\documentclass{amsart}
\usepackage{graphics,float}
\usepackage{amsmath,amsfonts,amssymb,amsthm}
\theoremstyle{plain}
\newtheorem{theorem}{Theorem}[section]
\newtheorem{conjecture}[theorem]{Conjecture}
\newtheorem{heuristic}[theorem]{Heuristic}
\newtheorem{guess}[theorem]{Guess}

\DeclareFontEncoding{OT2}{}{} 
  \newcommand{\textcyr}[1]{%
    {\fontencoding{OT2}\fontfamily{wncyr}\fontseries{m}\fontshape{n}%
     \selectfont #1}}
\newcommand{\Sha}{{\mbox{\textcyr{Sh}}}}

\newcommand{\Beta}{{B}}
\def\SO{{\rm SO}}
\def\agm{{\rm agm}}

\title{Some heuristics about elliptic curves}
\author{Mark Watkins}
\begin{document}
\begin{abstract}
We give some heuristics for counting elliptic curves with certain
properties. In particular, we re-derive the Brumer-McGuinness
heuristic for the number of curves with positive/negative discriminant
up to~$X$, which is an application of lattice-point counting.
We then introduce heuristics (with refinements from random matrix theory)
that allow us to predict how often we expect an elliptic curve~$E$
with even parity to have $L(E,1)=0$. We find that we expect there to be
about $c_1X^{19/24}(\log X)^{3/8}$ curves with $|\Delta|<X$
with even parity and positive (analytic) rank;
since Brumer and McGuinness predict $cX^{5/6}$ total curves,
this implies that asymptotically almost all even parity curves have rank~0.
We then derive similar estimates for ordering by conductor,
and conclude by giving various data regarding our heuristics
and related questions.
\end{abstract}
\maketitle

\vphantom 0\vskip-48pt\vphantom 0
\section{Introduction}
We give some heuristics for counting elliptic curves with certain
properties. In particular, we re-derive the Brumer-McGuinness
heuristic for the number of curves with positive/negative discriminant
up to~$X$, which is an application of lattice-point counting.
We then introduce heuristics (with refinements from random matrix theory)
that allow us to predict how often we expect an elliptic curve~$E$
with even parity to have $L(E,1)=0$. It turns out that we roughly expect
that a curve with even parity has $L(E,1)=0$ with probability proportional
to the square root of its real period, and, since the real period
is very roughly~$1/\Delta^{1/12}$, this leads us to the
prediction that almost all curves with even parity should have~$L(E,1)\neq 0$.
By the conjecture of Birch and Swinnerton-Dyer, this says that
almost all such curves have rank~0.

We then make similar heuristics when enumerating by conductor.
The first task here is simply to count curves with conductor up to~$X$,
and for this we use heuristics involving how often large powers of
primes divide the discriminant.
Upon making this estimate, we are then able to imitate
the argument we made previously,
and thus derive an asymptotic for the number of curves
with even parity and $L(E,1)=0$ under the ordering by conductor.
We again get the heuristic that almost all curves
with even parity should have~$L(E,1)\neq 0$.

We then make a few remarks regarding how often curves should
have nontrivial isogenies and/or torsion under different orderings,
and then present some data regarding average ranks. We conclude by
giving data for Mordell-Weil lattice distribution for rank 2 curves,
and speculating about symmetric power $L$-functions.

\section{The Brumer-McGuinness Heuristic}

First we re-derive the Brumer-McGuinness
heuristic~\cite{brumer-mcguinness} for the number of elliptic
curves whose absolute discriminant is less than a given bound~$X$;
the technique here is essentially lattice-point counting, and we
derive our estimates via the assumption that these counts are
well-approximated by the areas of the given regions.

\begin{conjecture}\label{conj:bmcg}[Brumer-McGuinness]
The number $A_\pm(X)$ of rational elliptic curves with a global
minimal model (including at $\infty$) and positive or negative
discriminant whose absolute value is less than~$X$
is asymptotically $A_\pm(X)\sim {\alpha_{\pm}\over\zeta(10)}X^{5/6}$, where
$\alpha_{\pm}={\sqrt 3\over 10}\int_{\pm1}^\infty {dx\over\sqrt{x^3\mp1}}$.
\end{conjecture}

As indicated by Brumer and McGuinness, the identity
$\alpha_{-}=\sqrt{3}\alpha_{+}$ was already known to Legendre,
and is related to complex multiplication.
These constants can be expressed in terms of Beta integrals
$B(u,v)=\int_0^1 x^{u-1}(1-x)^{v-1}\, dx=
{\Gamma(u)\Gamma(v)\over \Gamma(u+v)}$ as
$\alpha_{+}={1\over 3}\Beta(1/2,1/6)$ and $\alpha_{-}=\Beta(1/2,1/3)$.

Recall that every rational elliptic curve has a unique integral
minimal model $y^2+a_1xy+a_3y=x^3+a_2x^2+a_4x+a_6$ with
$a_1,a_3\in\{0,1\}$ and $|a_2|\le 1$.
Fix one of the 12 choices of~$(a_1,a_2,a_3)$. Since these are all
bounded by~1 the discriminant is thus approximately~$-64a_4^3-432a_6^2$.
So we wish to count the number
of $(a_4,a_6)$-lattice-points with~$|64a_4^3+432a_6^2|\le X$,
noting that Brumer and McGuinness divide the curves according to the
sign of the discriminant. The lattice-point count
for $a_1=a_2=a_3=0$ is given by
$$\mathop{\sum\sum}_{0<-64a_4^3-432a_6^2<X} 1
\quad+\mathop{\sum\sum}_{-X<-64a_4^3-432a_6^2<0} 1.$$
We estimate this lattice-point count by the integral
$\mathop{\int\int}_{U} du_4 du_6$ for the region $U$ given by
$|64u_4^3+432u_6^2|<X$. After splitting into two parts based
upon the sign of the discriminant and performing the $u_4$-integration,
we get
\begin{align*}
{2\over (64)^{1/3}}
\int_0^\infty \Bigl[(-432u_6^2)^{1/3}&-(-X-432u_6^2)^{1/3}\Bigr]\, du_6+\\
&+{2\over (64)^{1/3}}
\int_0^\infty \Bigl[(X-432u_6^2)^{1/3}-(-432u_6^2)^{1/3}\Bigr]\, du_6,
\end{align*}
where the factor of $2$ comes from the sign of~$u_6$.
Changing variables $u_6=w\sqrt{X/432}$ and multiplying by~12
for the number of $(a_1,a_2,a_3)$-choices we get
\begin{align*}
{24\over (64)^{1/3}}{X^{5/6}\over \sqrt{432}}
\int_0^\infty \Bigl[(w^2+1)^{1/3}&-(w^2)^{1/3}\Bigr]\, dw+\\
&+{24\over (64)^{1/3}}{X^{5/6}\over \sqrt{432}}
\int_0^\infty \Bigl[(w^2)^{1/3}-(w^2-1)^{1/3}\Bigr]\, dw.
\end{align*}
These integrals are probably known, but I am unable to find a reference.
The first integral simplifies\footnote{
 As N.~D.~Elkies indicated to us, we can write
 $I(a)=\int_0^\infty \bigl[(t^2+a)^{1/3}-(t^2)^{1/3}\bigr]\, dt$,
 differentiate under the integral sign, then substitute~$t^2+a=ax^3$,
 and finally re-integrate to obtain~$I(1)$.}
to ${3\over 5}\int_{1}^\infty {dx\over\sqrt{x^3-1}}={1\over 5}\Beta(1/2,1/6)$,
while the second becomes
${3\over 5}\int_{-1}^\infty {dx\over\sqrt{x^3+1}}={3\over 5}\Beta(1/2,1/3)$.
This counts all models of curves; if we eliminate non-minimal models,
for which we have $p^4|c_4$ and~$p^{12}|\Delta$ for some prime~$p$,
we expect to accrue an extra factor\footnote{
 Note that some choices of $(a_1,a_2,a_3)$ necessarily
 have odd discriminant, but the other choices compensate
 to give the proper Euler factors at 2 (and 3).}
of~$\zeta(10)$. From this, we get the conjecture of Brumer and McGuinness
as stated above.

\section{Counting curves of even parity whose central $L$-value vanishes.}

See \cite[\S 15-16]{silverman} for definitions of the conductor $N$
and $L$-function $L(E,s)$ of an elliptic curve~$E$.
Since rational elliptic curves are modular, we have that the
completed $L$-function $\Lambda(E,s)=\Gamma(s)(\sqrt N/2\pi)^sL(E,s)$
satisfies $\Lambda(E,s)=\pm\Lambda(E,2-s)$. When the plus sign occurs,
we say that $E$ has even parity.

We now try to count elliptic curves $E$ with
even parity for which~$L(E,1)=0$.
Throughout this section, $E$ shall be a curve with
even parity, and we shall order curves by discriminant.
Via the conjectural Parity Principle, we expect that,
under any reasonable ordering, half of the elliptic curves
should have even parity; in particular, we predict that there
are asymptotically $A_{\pm}(X)/2$ curves with even parity and
positive/negative discriminant up to~$X$.

Our main tool shall be random matrix theory, which
gives a heuristic for predicting how often $L(E,1)$ is small.
We could alternatively derive a cruder heuristic by assuming the
the order of the Shafarevich-Tate group is a random square
integer in a given interval, but random matrix theory has the
advantage of being able to predict a more explicit asymptotic.
Our principal heuristic is the following:

\begin{heuristic}\label{conj:rank2disc}
The number $R(X)$ of rational elliptic curves $E$ with
even parity and $L(E,1)=0$ and
absolute discriminant less than~$X$
is given asymptotically by $R(X)\sim cX^{19/24}(\log X)^{3/8}$
for some computable constant~$c>0$.
\end{heuristic}

In particular, note that we get the prediction
that almost all curves with
even parity have $L(E,1)\neq 0$ under this ordering.

\subsection{Random matrix theory}
Originally developed in mathematical statistics by Wishart \cite{wishart}
in the 1920s and then in mathematical physics
(especially the spectra of highly excited nuclei)
by Wigner \cite{wigner}, Dyson, Mehta, and others
(particularly~\cite{marcenko-pastur}), random matrix theory \cite{mehta}
has now found some applications in number theory,
the earliest being the oft-told story of
Dyson's remark to Montgomery regarding the pair-correlation of zeros of
the Riemann $\zeta$-function.
Based on substantial numerical evidence, random matrix theory appears
to give reasonable models for the distribution of $L$-values in families,
though the issue of what constitutes a proper family is a delicate one
(see particularly \cite[\S 3]{CFKRS}, where the notion of family comes from
the ability to calculate moments of $L$-functions rather than
from algebraic geometry).

The family of quadratic twists of a given elliptic curve
$E:y^2=x^3+Ax+B$ is given by $E_d:y^2=x^3+Ad^2x+Bd^3$ for squarefree~$d$.
The work (most significantly a monodromy computation)
of Katz and Sarnak \cite{katz-sarnak} regarding families of
curves over function fields implies that when we restrict to quadratic twists
with even parity, we should expect that the $L$-functions are modelled
by random matrices with even orthogonal symmetry.
Though we have no function field analogue in our case, we brazenly
assume (largely from looking at the sign in the functional equation)
that the symmetry type is again orthogonal with even parity.
What this means is that we want to model properties of the $L$-function 
via random matrices taken from $\SO(2M)$ with respect to Haar measure.
Here we wish the mean density of zeros
of the \hbox{$L$-functions} to match the mean density of eigenvalues
of our matrices, and so, as in~\cite{keating-snaith},
we should take $2M\approx 2\log N$. We suspect that the $L$-value
distribution is approximately given by the distribution of the evaluations
at~$1$ of the characteristic polynomials of our random matrices.
In the large, this distribution is determined entirely by
the symmetry type, while finer considerations are distinguished
via arithmetic considerations.

With this assumption, via the moment conjectures of \cite{keating-snaith}
and then using Mellin inversion, as $t\rightarrow 0$ we have
(see (21) of \cite{random-matrix-theory}) that
\begin{equation}\label{RMTprob}
{\rm Prob}[L(E,1)\le t]\sim \alpha(E) t^{1/2}M^{3/8}.
\end{equation}
This heuristic is stated for fixed $M\approx\log N$,
but we shall also allow $M\rightarrow\infty$.
It is not easy to understand this probability, as both the constant
$\alpha(E)$ and the matrix-size $M$ depend on~$E$.
We can take curves with $e^M\le N\le e^{M+1}$ to mollify the
impact of the conductor, but in order to average over a set of
curves, we need to understand how $\alpha(E)$ varies.
One idea is that $\alpha(E)$ separates into two parts, one of which
depends on local structure (Frobenius traces) of the curve,
and the other of which depends only upon the size of the conductor~$N$.
Letting $G$ be the Barnes $G$-function (such that $G(z+1)=\Gamma(z)G(z)$
with $G(1)=1$) and $M=\lfloor\log N\rfloor$ we have that
$$\alpha(E)=\alpha_R(M)\cdot \alpha_A(E)$$
$\text{with}\>\>\alpha_R(M)
\rightarrow\hat\alpha_R=2^{1/8}G(1/2)\pi^{-1/4}\>\>
\text{as}\>\> M\rightarrow\infty$
and
\begin{equation}\label{eqn:Fp}
\alpha_A(E)=\prod_p F(p)=
\prod_p \biggl(1-{1\over p}\biggr)^{3/8} \biggl({p\over p+1}\biggr)
\biggl({1\over p}+{L_p(1/p)^{-1/2}\over 2}+{L_p(-1/p)^{-1/2}\over 2}\biggr)
\end{equation}
where $L_p(X)=(1-a_pX+pX^2)^{-1}$ when $p\nmid\Delta$
and $L_p(X)=(1-a_pX)^{-1}$ otherwise;
see (10) of \cite{random-matrix-theory} evaluated at~$k=-1/2$,
though that equation is wrong at primes that
divide the discriminant --- see (20) of~\cite{CPRW},
where $Q$ should be taken to be~1.
Note that the Sato-Tate conjecture \cite{tate} implies that
$a_p^2$ is $p$ on average, and this implies that the above
Euler product converges.

\subsection{Discretisation of the $L$-value distribution}
For precise definitions of the Tamagawa numbers, torsion group,
periods, and Shafarevich-Tate group, see~\cite{silverman},
though below we give a brief description of some of these.
We let $\tau_p(E)$ be the Tamagawa number of~$E$ at the
(possibly infinite) prime~$p$, and write $\tau(E)=\prod_p \tau_p(E)$
for the Tamagawa product and $T(E)$ for the size of the torsion group.
We also write $\Omega_{\rm re}(E)$ for the real
period and $\Sha_{\rm an}(E)$ for the size of the Shafarevich-Tate group
when $L(E,1)\neq 0$, with $\Sha_{\rm an}(E)=0$ when $L(E,1)=0$.

We wish to assert that sufficiently small values of~$L(E,1)$
actually correspond to~$L(E,1)=0$.
We do this via the conjectural formula of Birch and
Swinnerton-Dyer~\cite{BSD}, which asserts that
$$L(E,1)=
\Omega_{\rm re}(E)\cdot {\tau(E)\over T(E)^2}\cdot \Sha_{\rm an}(E).$$
Our discretisation\footnote
{The precision of this discretisation might be the most-debatable
 methodology we use. Indeed, we are essentially taking a ``sharp cutoff'',
 while it might be better to have a more smooth transition function.
 For this reason, we do not specify the leading
 constant in our final heuristic.}
will be that
$$L(E,1)<\Omega_{\rm re}(E)\cdot {\tau(E)\over T(E)^2}
\quad\text{implies}\quad L(E,1)=0.$$
Note that we are only using that $\Sha_{\rm an}$ takes on integral
values, and do not use the (conjectural) fact that it is square.

Using \eqref{RMTprob}, we estimate the number of curves
with positive (for simplicity) discriminant less than~$X$
and even parity and $L(E,1)=0$ via the lattice-point sum
$$W(X)=\mathop{\sum\sum}_{\text{$c_4,c_6$ minimal}\atop 0<c_4^3-c_6^2<1728X}
\alpha_R(M)\alpha_A(E)
\cdot{\sqrt{\Omega_{\rm re}(E)\tau(E)\over T(E)^2}}\cdot M^{3/8}.$$
We need to introduce congruence conditions on $c_4$ and $c_6$
to make sure that they correspond to a minimal model of an elliptic curve.
The paper \cite{stein-watkins} uses the work of Connell \cite{connell}
in a different context to get that
there are 288 classes of $(c_4\>\text{mod}\>576,c_6\>\text{mod}\>1728)$
that can give minimal models, and so we get a factor of~$288/(576\cdot 1728)$,
assuming that each congruence class has the same impact
on all the entities in the sum. Indeed, this independence
(on average) of various quantities with respect
to $c_4$ and $c_6$ is critical in our estimation of~$W(X)$.
There is also the question of non-minimal models,\footnote
{At $p=2,3$, non-minimality occurs when $c_4/p^4$ and $c_6/p^6$
 satisfy the congruences.}
from which we get a factor of~$1/\zeta(10)$.

\begin{guess}\label{GUESS}
The lattice-point sum $W(X)$ can be approximated as $X\rightarrow\infty$ by
$$\hat W(X)={288\over (576\cdot 1728)}{1\over\zeta(10)}\cdot
\hat\alpha_R\bar\alpha_A\beta_\tau\cdot
\mathop{\int\int}_{1\le{u_4^3-u_6^2\over 1728}<X}
\Omega_{\rm re}(E)^{1/2}\cdot (\log\Delta)^{3/8} du_4du_6.$$
Here $\hat\alpha_R$ is the limit $2^{1/8}G(1/2)\pi^{-1/4}$ of $\alpha_R(M)$
as $M\rightarrow\infty$, while $\bar\alpha_A$ is a suitable average
of the arithmetic factors $\alpha_A(E)$, and $\beta_\tau$ is the
average of the square root of the Tamagawa product. We have also approximated
$\log N\approx\log\Delta$ and assumed the torsion is trivial;
below we will give these heuristic justification~(on average).
Note that everything left in the integral is
a smooth function of $u_4$ and $u_6$.
\end{guess}

We shall first evaluate the integral in $\hat W(X)$ given these suppositions,
and then try to justify the various assumptions that are inherent in
this guess.\footnote
{Note that our methods do not readily generalise to higher rank, as there is
 no apparent way to model the heights of points (and thus the regulator).}

\subsection{Evaluation of the integral}
Write $E$ as $y^2=4x^3-(c_4/12)x-c_6/216$,
and put $e_1>e_2>e_3$ for the roots of the cubic polynomial on the right side.
We have
$$1/\Omega_{\rm re}=
{\rm agm}\bigl(\sqrt{e_1-e_2},\sqrt{e_1-e_3}\bigr)/\pi.$$
We also have that $(e_1-e_2)(e_1-e_3)(e_2-e_3)=\sqrt{\Delta/16}$
from the formula for the discriminant.
We next write $e_1-e_2=\Delta^{1/6}\lambda$ and $e_2-e_3=\Delta^{1/6}\mu$
so that we have $\mu\lambda(\lambda+\mu)=1/4$,
while $e_1={\Delta^{1/6}\over 3}(\mu+2\lambda)$,
$e_2={\Delta^{1/6}\over 3}(\mu-\lambda)$,
and $e_3=-{\Delta^{1/6}\over 3}(2\mu+\lambda)$.
Thus we get
$$-c_6/864=-e_1e_2e_3=
{\Delta^{1/2}\over 27}(\mu+2\lambda)(\mu-\lambda)(2\mu+\lambda)$$
and
$$-c_4/48=e_1e_2+e_1e_3+e_2e_3=
-{\Delta^{1/3}\over 3}(\mu^2+\lambda\mu+\lambda^2).$$
Changing variables in the $\hat W$-integral
gives a Jacobian of $432/\Delta^{1/6}\sqrt{\mu^4+\mu}$
so that
$$\hat W(X)=\tilde c\int_1^X\int_0^\infty {(\log \Delta)^{3/8}\over
\sqrt{\Delta^{1/12}\,{\rm agm}(\sqrt\lambda,\sqrt{\lambda+\mu})}}
{d\mu\,d\Delta\over\Delta^{1/6}\sqrt{\mu^4+\mu}},$$
where $\lambda=(\sqrt{\mu^4+\mu}-\mu^2)/2\mu$.
Thus the variables are nicely separated, and since the $\mu$-integral
converges, we do indeed get $\hat W(X)\sim cX^{19/24}(\log X)^{3/8}$.
A similar argument can be given for curves with negative discriminant.
This concludes our derivation of Heuristic~\ref{conj:rank2disc},
and now we turn to giving some reasons for our expectation that
the arithmetic factors can be mollified by taking their averages.

\subsection{Expectations for arithmetic factors on average}
In the next section we shall explain (among other things)
why we expect that $\log N\approx \log\Delta$ for almost all curves,
and in section~\ref{section:torsion},
we shall recall the classical parametrisations of $X_1(N)$
due to Fricke to indicate why we expect the torsion size
is 1 on average. Here we show how to compute the various averages
(with respect to ordering by discriminant)
of the square root of the Tamagawa product
and the arithmetic factors~$\alpha_A(E)$.

For both heuristics, we shall make the assumption that curves satisfying
the discriminant bound $|\Delta|\le X$ behave essentially the same as
those that satisfy $|c_4|\le X^{1/3}$ and $|c_6|\le X^{1/2}$.
That is, we approximate our region by a big box.
We write $D$ for the absolute value of~$\Delta$.
First we consider the Tamagawa product.

We wish to know how often a prime divides the discriminant to a high power.
Fix a prime~$p\ge 5$ with $p$ a lot smaller than~$X^{1/3}$.
We can estimate the probability that~$p^k|\Delta$
by considering all $p^{2k}$ choices of $c_4$ and $c_6$ modulo~$p^k$,
that is, by counting the number of solutions $S(p^k)$
to $c_4^3-c_6^2=1728\Delta\equiv 0$ (mod~$p^k$).
This auxiliary curve $c_4^3=c_6^2$ is singular at $(0,0)$ over~${\bf F}_p$,
and has $(p-1)$ non-singular ${\bf F}_p$-solutions
which lift to $p^{k-1}(p-1)$ points modulo~$p^k$.

For $p^k$ sufficiently small, our $(c_4,c_6)$-region is so large that we can 
show that the probability that $p^k|\Delta$ is~$S(p^k)/p^{2k}$.
We assume that big primes act (on average) in the same manner,
while a similar heuristic can be given for~$p=2,3$.
Curves with $p^4|c_4$ and $p^6|c_6$ will not be given by their
minimal model; indeed, we want to exclude these curves, and
thus will multiply our probabilities by $\kappa_p=(1-1/p^{10})^{-1}$
to make them conditional on this criterion. For instance, the above
counting of points says that there is a probability of $(p^2-p)/p^2$
that~$p\nmid D$, and so upon conditioning upon minimal models we get
$\kappa_p(1-1/p)$ for this probability.

What is the probability $P_m(p,k)$ that a curve given by a minimal model
has multiplicative reduction at $p\ge 5$ and $p^k\|D$ for some~$k>0$?
In terms of Kodaira symbols, this is the case of I$_k$.
For multiplicative reduction we need that $p\nmid c_4,c_6$.
These events are independent and each has a probability $(1-1/p)$
of occurring. Upon assuming these conditions and working modulo~$p^k$,
there are $(p^k-p^{k-1})$ such choices for each,
and of the resulting $(c_4,c_6)$
pairs we noted above that $p^{k-1}(p-1)$ of them have~$p^k|D$.
So, given a curve with~$p\nmid c_4,c_6$,
we have a probability of $1/p^{k-1}(p-1)$ that~$p^k|D$,
which gives $1/p^k$ for the probability that~$p^k\|D$.
In symbols, we have that (for~$p\ge 5$ and $k\ge 1$)
$${\rm Prob}\Bigl[p^k\|(c_4^3-c_6^2) \Bigm| p\nmid c_4,c_6\Bigr]=1/p^k.$$
Including the conditional probability for minimal models, we get
\begin{equation}\label{eqn:probm}
P_m(p,k)=(1-1/p^{10})^{-1}(1-1/p)^2/p^k\quad\text{(for $p\ge 5$ and $k\ge 1$).}
\end{equation}
Note that summing this over $k\ge 1$ gives $\kappa_p(1-1/p)/p$ for the
probability for an elliptic curve to have multiplicative reduction at~$p$.

What is the probability $P_a(p,k)$ that a curve given by a minimal
model has additive reduction at $p\ge 5$ and $p^k\|D$ for some~$k>0$?
We shall temporarily ignore the factor
of $\kappa_p=(1-1/p^{10})^{-1}$ from non-minimal models 
and include it at the end.
We must have that $p|c_4,c_6$, and thus get that~$k\ge 2$.
For $k=2,3,4$, which correspond to Kodaira symbols II, III, and IV
respectively, the computation is not too bad:
we get that $p^2\|D$ exactly when $p|c_4$ and~$p\|c_6$,
so that the probability is $(1/p)\cdot (1-1/p)/p=(1-1/p)/p^2$;
for $p^3\|D$ we need $p\|c_4$ and $p^2|c_6$
and thus get $(1-1/p)/p\cdot(1/p^2)=(1-1/p)/p^3$ for the probability;
and for $p^4\|D$ we need $p^2|c_4$ and $p^2\|c_6$ and so get
$(1/p^2)\cdot(1-1/p)/p^2=(1-1/p)/p^4$ for the probability.
Note that the case $k=5$ cannot occur.
Thus we have (for~$p\ge 5$) the formula
$P_a(p,k)=(1-1/p^{10})^{-1}(1-1/p)/p^k$ for $k=2,3,4$.

More complications occur for $k\ge 6$, where now we split into two
cases depending upon whether additive reduction persists upon taking
the quadratic twist by~$p$. This occurs when $p^3|c_4$ and $p^4|c_6$,
and we denote by $P_a^n(p,k)$ the probability that $p^k\|D$ in this subcase.
Just as above, we get that
$P_a^n(p,k)=(1-1/p^{10})^{-1}(1-1/p)/p^{k-1}$ for $k=8,9,10$.
These are respectively the cases of Kodaira symbols IV$^\star$,
III$^\star$, and~II$^\star$.
For $k=11$ we have~$P_a^n(p,k)=0$, while for $k\ge 12$ our
condition of minimality implies that we should take~$P_a^n(p,k)=0$.

We denote by $P_a^t(p,k)$ the probability that $p^6|D$ with either
$p^2\|c_4$ or $p^3\|c_6$. First we consider curves for which $p^7|D$,
and these have multiplicative reduction at~$p$ upon twisting.
In particular, these curves have $p^2\|c_4$ and $p^3\|c_6$,
and the probability of this is $(1-1/p)/p^2\cdot (1-1/p)/p^3$.
Consider $k\ge 7$. We then take $c_4/p^2$ and $c_6/p^3$ both modulo $p^{k-6}$,
and get that $p^{k-6}\|(D/p^6)$ with probability $1/p^{k-6}$
in analogy with the above.
So we get that $P_a^t(p,k)=(1-1/p^{10})^{-1}(1-1/p)^2/p^{k-1}$ for~$k\ge 7$.
This corresponds to the case of I$_{k-6}^\star$.

Finally, for $p^6\|D$ (which is the case I$_0^\star$)
we get a probability of $(1/p^2)\cdot(1/p^3)$ for the chance
that $p^2|c_4$ and $p^3|c_6$, and (since there are $p$ points mod~$p$
on the auxiliary curve $(c_4/p^2)^3\equiv(c_6/p^3)^2 \pmod{p}$)
a conditional probability of $(p^2-p)/p^2$ that~$p^6\|D$.
So we get that $P_a^t(p,6)=(1-1/p^{10})^{-1}(1-1/p)/p^5$.

We now impose our current notation on the previous paragraphs,
and naturally let $P_a^t(p,k)=0$ and $P_a^n(p,k)=P_a(p,k)$ for~$k\le 5$.
Our final result is that
\begin{equation}\label{eqn:proba}
P_a^n(p,k)=
\begin{cases}
(1-1/p^{10})^{-1}(1-1/p)/p^k & \quad k=2,3,4\\
(1-1/p^{10})^{-1}(1-1/p)/p^{k-1} & \quad k=8,9,10
\end{cases}
\end{equation}
\begin{equation}\label{eqn:probd}
P_a^t(p,k)=
\begin{cases}
(1-1/p^{10})^{-1}(1-1/p)/p^5 & \quad k=6 \\
(1-1/p^{10})^{-1}(1-1/p)^2/p^{k-1} & \quad k\ge 7
\end{cases}
\end{equation}
with $P_a^n(p,k)$ and $P_a^t(p,k)$ equal to zero for other~$k$.
We conclude by defining $P_0(p,k)$ to be zero for $k>0$
and to be the probability $(1-1/p^{10})^{-1}(1-1/p)$ that $p\nmid D$ for~$k=0$.
We can easily check that we really do have the required probability relation
$\sum_{k=0}^\infty \bigl[P_m(p,k)+P_a^n(p,k)+P_a^t(p,k)+P_0(p,k)\bigr]=1$,
as: the cases of multiplicative reduction give $\kappa_p(1-1/p)/p$;
the cases of Kodaira symbols II, III, and~IV give $\kappa_p(1/p^2-1/p^5)$;
the cases of Kodaira symbols IV$^\star$, III$^\star$, and~II$^\star$
give $\kappa_p(1/p^7-1/p^{10})$; the cases of I$^\star_k$ summed
for~$k\ge 1$ give $\kappa_p(1-1/p)/p^6$; the case of I$^\star_0$
gives $\kappa_p(1-1/p)/p^5$; and the sum of these with
$P_0(p,0)=\kappa_p(1-1/p)$ does indeed give us~1.
We could do a similar (more tedious) analysis for $p=2,3$,
but this would obscure our argument.

Given a curve of discriminant~$D$,
we can now compute the expectation for its Tamagawa number.
We consider primes~$p|D$ with $p\ge 5$,
and compute the local Tamagawa number~$t(p)$.
When $E$ has multiplicative reduction at~$p$ and~$p^k\|D$,
then $t(p)=k$ if $-c_6$ is square mod~$p$,
and else $t(p)=1,2$ depending upon whether $k$ is odd or even.
So the average of $\sqrt{t(p)}$ for this case is
$\epsilon_m(k)={1\over 2}(1+\sqrt k),{1\over 2}(\sqrt 2+\sqrt k)$
for $k$ odd/even respectively.

When $E$ has potentially multiplicative reduction at~$p$ with~$p^k\|D$,
for $k$ odd we have $t(p)=4,2$ depending on whether
$(c_6/p^3)\cdot(\Delta/p^k)$ is square mod~$p$,
and for $k$ even we have $t(p)=4,2$ depending on
whether $\Delta/p^k$ is square mod~$p$.
In both cases the average of $\sqrt{t(p)}$ is~${1\over 2}(\sqrt 2+\sqrt 4)$.
In the case of I$_0^\star$ reduction where we have~$p^6\|D$,
we have that $t(p)=1,2,4$ corresponding to whether the cubic
$x^3-(27c_4/p^2)x-(54c_6/p^3)$ has $0,1,3$ roots modulo~$p$.
So the average of $\sqrt{t(p)}$ is
$${\sqrt 1\bigl((p-1)(p+1)/3\bigr)+\sqrt 2\bigl(p(p-1)/2\bigr)+
\sqrt 4\bigl((p-1)(p-2)/6\bigr)\over
\bigl((p-1)(p+1)/3\bigr)+\bigl(p(p-1)/2\bigr)+
\bigl((p-1)(p-2)/6\bigr)}={2\over 3}+{\sqrt 2\over 2}-{1\over 3p}.$$
in this case.

For the remaining cases,
when $p^2\|D$ or $p^{10}\|D$ we have~$t(p)=1$,
while when $p^3\|D$ or $p^9\|D$ we have~$t(p)=2$.
Finally, when $p^4\|D$ we have $t(p)=3,1$ depending
on whether $-6c_6/p^2$ is square mod~$p$,
and similarly when $p^8\|D$ we have $t(p)=3,1$ depending
on whether $-6c_6/p^4$ is square mod~$p$, so that the average
of $\sqrt{t(p)}$ in both cases is~${1\over 2}(1+\sqrt 3)$.
We get that $\epsilon_a^n(k)=1,\sqrt 2,{1\over 2}(1+\sqrt 3),
{1\over 2}(1+\sqrt 3),\sqrt 2,1$ for $k=2,3,4,8,9,10$, while
\begin{equation}\label{eqn:epsilon}
\epsilon_m(k)=
\begin{cases}{1\over 2}(1+\sqrt k),&\text{$k$ odd}\\
{1\over 2}(\sqrt 2+\sqrt k),&\text{$k$ even}\end{cases}
\quad\text{and}\quad
\epsilon_a^t(p,k)=\begin{cases}
{2\over 3}+{\sqrt 2\over 2}-{1\over 3p},&k=6\\
{1\over 2}(\sqrt 2+\sqrt 4),&k\ge 7\end{cases}
\end{equation}
with $\epsilon_a^n(k)$ and $\epsilon_a^t(p,k)$ equal to zero for other~$k$.

We define the expected square root of the Tamagawa number $K(p)$ at~$p$ by
\begin{equation}\label{eqn:tama}
K(p)=\sum_{k=0}^\infty
\bigl[\epsilon_m(k) P_m(p,k)+\epsilon_a^n(k) P_a^n(p,k)+
\epsilon_a^t(p,k) P_a^t(p,k)+P_0(p,k)\bigr]
\end{equation}
and the expected global\footnote
{Note that the Tamagawa number at infinity is 1 when $E$
 has negative discriminant and else is 2, the former
 occurring approximately $\sqrt{3}/(1+\sqrt{3})\approx 63.4\%$ of the time.}
Tamagawa number to be~$\beta_\tau=\prod_p K(p)$.
The convergence of this product follows from an analysis of the
dominant $k=0,1,2$ terms of~\eqref{eqn:tama}, which gives 
a behaviour of~$1+O(1/p^2)$.
So we get that the Tamagawa product is a constant on average,
which we do not bother to compute explicitly (we would need to
consider $p=2,3$ more carefully to get a precise value).

To compute the average value of $\alpha_A(E)=\prod_p F(p)$
in~\eqref{eqn:Fp} we similarly assume\footnote
{This argumentative technique can also be used to bolster our assumption that
 using Connell's conditions should be independent of other considerations.}
that each prime acts independently; we then compute the
average value for each prime by calculating the distribution of $F(p)$
when considering all the curves modulo~$p$
(including those with singular reduction,
and again making the slight adjustment for non-minimal models).
This gives some constant for the average $\bar\alpha_A$ of~$\alpha_A(E)$,
which we do not compute explicitly.
Note that $\prod_p F(p)$ converges if we assume the Sato-Tate
conjecture \cite{tate} since in this case we have that
$a_p^2$ is $p$ on average.

\section{Relation between conductor and discriminant}
We now give heuristics for how often we expect the ratio
between the absolute discriminant and the conductor to be large.
The main heuristic we derive is:

\begin{heuristic}\label{conj:conductor}
The number $B(X)$ of rational elliptic curves whose conductor
is less than~$X$ satisfies $B(X)\sim cX^{5/6}$ for an explicit
constant~$c>0$.
\end{heuristic}

To derive this heuristic, we estimate the proportion of curves with a given
ratio of (absolute) discriminant to conductor. Since the conductor
is often the squarefree kernel of the discriminant, by way of
explanation we first consider the behaviour of $f(n)=n/{\rm sqfree}(n)$.
The probability that $f(n)=1$ is given by the probability that $n$
is squarefree, which is classically known to be $1/\zeta(2)=6/\pi^2$.
Given a prime power~$p^m$, to have $f(n)=p^m$ says that $n=p^{m+1}u$
where $u$ is squarefree and coprime to~$p$. The probability that
$p^{m+1}\|n$ is $(1-1/p)/p^{m+1}$, and given this, the conditional
probability that $\bigl(n/p^{m+1}\bigr)$
is squarefree is $(6/\pi^2)\cdot(1-1/p^2)^{-1}$.
Extending this multiplicatively beyond prime powers, we get that
$${\rm Prob}\bigl[n/{\rm sqfree}(n)=q\bigr]=
{6\over\pi^2}\prod_{p^m\|q} {1/p^{(m+1)}\over (1+1/p)}=
{6\over\pi^2}{1\over q}\prod_{p|q} {1\over p+1}.$$
In particular, the average of $f(n)^\gamma$ exists for~$\gamma<1$;
in our elliptic curve analogue, we will require such
an average for~$\gamma=5/6$. We note that it appears to be
an interesting open question to prove an asymptotic
for~$\sum\limits_{n\le X} n/{\rm sqfree}(n)$.

\subsection{Derivation of the heuristic}
We keep the notation $D=|\Delta|$ and wish to compute
the probability that $D/N=q$ for a fixed positive integer~$q$.
For a prime power~$p^v$ with $p\ge 5$, the probability that $p^v\|(D/N)$
is given by: the probability that $E$ has multiplicative reduction
at~$p$ and~$p^{v+1}\|D$, that is~$P_m(p,v+1)$; plus the probability
that $E$ has additive reduction at~$p$ and~$p^{v+2}\|D$, that is~$P_a(p,v+2)$;
and the contribution from~$P_0(p,v)$, which is zero for $v>0$
and for $v=0$ is the probability that $p$ does not divide~$D$.
So, writing $v=v_p(q)$, we get that
(with a similar modified formula for~$p=2,3$)
\begin{equation}\label{eqn:prob}
{\rm Prob}\bigl[D/N=q\bigr]=
\prod_p E_p(v_p(q))=\prod_p\bigl[P_m(p,1+v)+P_a(p,2+v)+P_0(p,v)\bigr].
\end{equation}
It should be emphasised that this probability is with respect to
(as in the previous section) the ordering of the curves by discriminant.
We have
\begin{equation}\label{eqn:sum}
\sum_{E: N_E\le X} 1
\approx\sum_{q=1}^\infty\sum_{E: D\le qX} {\rm Prob}\bigl[D/N=q\bigr]
\sim\sum_{q=1}^\infty\alpha (qX)^{5/6}\cdot{\rm Prob}\bigl[D/N=q\bigr],
\end{equation}
where $\alpha=\alpha_++\alpha_-$
from the Brumer-McGuinness heuristic~\ref{conj:bmcg}.
If this last sum converges, then we get Heuristic~\ref{conj:conductor}.

To show the last sum in~\eqref{eqn:sum} does indeed converge,
we upper-bound the probability in~\eqref{eqn:prob}.
We have that $P_m(p,v+1)\le 1/p^{v+1}$ and~$P_a(p,v+2)\le 2/p^{v+1}$,
which implies
$$\hat f(q)={\rm Prob}\bigl[D/N=q\bigr]=
\prod_p E_p(v_p(q))\le {1\over q}\prod_{p|q} {3\over p}.$$
We then estimate
$$\sum_{q=1}^\infty q^{5/6}\hat f(q)\le
\sum_{q=1}^\infty {1\over q^{1/6}}\prod_{p|q} {3\over p}=
\prod_p\biggl(1+\sum_{l=1}^\infty {3/p\over (p^l)^{1/6}}\biggr)
\le\prod_p\biggl(1+{3/p\over p^{1/6}-1}\biggr),$$
and the last product is convergent upon comparison to~$\zeta(7/6)^3$.
Thus we shown that the last sum in~\eqref{eqn:sum} converges,
so that Heuristic~\ref{conj:conductor} follows.

We can note that Fouvry, Nair, and Tenenbaum \cite{FNT}
have shown that the number of minimal models with $D\le X$
is at least $cX^{5/6}$, and that the number of curves with
$D\le X$ with Szpiro ratio ${\log D\over\log N}\ge \kappa$
is no more than $c_\epsilon X^{1/\kappa+\epsilon}$ for every~$\epsilon>0$.

\subsection{Dependence of $D/N$ and the Tamagawa product}
We expect that $D/N$ should be independent of the real period,
but the Tamagawa product and $D/N$ should be somewhat related.\footnote
{The size of the torsion subgroup should also be related to~$D/N$,
 but in the next section we argue that curves with nontrivial
 torsion are sufficiently sparse so as to be ignored.}
We compute the expected square root of the Tamagawa product when~$D/N=q$.
As with~\eqref{eqn:prob} and using the $\epsilon$ defined
in~\eqref{eqn:epsilon}, we find that this is given by
$$\eta(q)=\prod_p
{\bigl[\epsilon_m(v_1)P_m(p,v_1)+\epsilon_a^n(v_2)P_a^n(p,v_2)+
\epsilon_a^t(p,v_2)P_a^t(p,v_2)+P_0(p,v)\bigr]
\over \bigl[P_m(p,v_1)+P_a(p,v_2)+P_0(p,v)\bigr]},$$
where $v_1=v+1$ and $v_2=v+2$ and $v=v_p(q)$.

\subsection{The comparison of $\log\Delta$ with $\log N$}
We now want to compare $\log\Delta$ with~$\log N$,
and explicate the replacement therein in Guess~\ref{GUESS}.
In order to bound the effect of curves with large~$D/N$, we note that
$${\rm Prob}\bigl[D/N\ge Y\bigr]=\sum_{q\ge Y} \hat f(q)
\le\sum_{q\ge Y} {1\over q}\prod_{p|q}{3\over p},$$
and use Rankin's trick, so that for any $0<\alpha<1$ we have
(using $p^\alpha-1\ge\alpha\log p$)
\begin{align*}
{\rm Prob}\bigl[D/N\ge Y\bigr]&
\le\sum_{q=1}^\infty
\biggl({q\over Y}\biggr)^{1-\alpha}\cdot {1\over q}\prod_{p|q}{3\over p}
={Y^\alpha\over Y}\prod_p\biggl(1+{3\over p^{1+\alpha}}+
{3\over p^{1+2\alpha}}+\cdots\biggr)\\
&={Y^\alpha\over Y}\prod_p\biggl(1+{3/p\over p^\alpha-1}\biggr)
\ll{Y^\alpha\over Y}\exp\biggl(\sum_p{\hat c/p\over \alpha\log p}\biggr)
\ll {e^{c\sqrt{\log Y}}\over Y}
\end{align*}
for some constants~$\hat c,c$, by taking $\alpha=1/\sqrt{\log Y}$
(this result is stronger than needed).

However, a more pedantic derivation of Guess \ref{GUESS} does not simply
allow replacing $\log N$ by~$\log\Delta$, but requires analysis
(assuming $\Omega_{\rm re}(E)$ to be independent of~$q$) of
$${\hat\alpha_R\hat\alpha_A\over 3456\,\zeta(10)}
\cdot\kern-10pt\mathop{\int\int}_{\sqrt X\le {u_4^3-u_6^2\over 1728}\le X}
\kern-10pt\Omega_{\rm re}(E)
\cdot\biggl[\sum_{q<\Delta} \eta(q)(\log \Delta/q)^{3/8}
\cdot{\rm Prob}\bigl[D/N=q\bigr]\biggr]\, du_4\,du_6.$$
The above estimate on the tail of the probability
and a simple bound on $\eta(q)$ in terms of the divisor function
shows that we can truncate the $q$-sum at~$Y$
with an error of~$O(1/Y^{8/9})$, and choosing (say) $Y=e^{\sqrt{\log X}}$
gives us that~$\log(\Delta/q)\sim\log\Delta$
(note that we restricted to~$\Delta>\sqrt X$).
So the bracketed term becomes the desired
$$\sum_{q<Y} (\log \Delta)^{3/8}
\eta(q)\cdot{\rm Prob}\bigl[D/N=q\bigr]\sim \beta_\tau(\log\Delta)^{3/8},$$
upon noting that the $q$-part of the sum converges to~$\beta_\tau$
as~$Y\rightarrow\infty$.

\subsection{Counting curves with vanishing L-value}
We now estimate the number of elliptic curves~$E$ with even parity
and $L(E,1)=0$ when ordered by conductor.

\begin{heuristic}\label{conj:rank2cond}
Let $\tilde R(X)$ be the number of elliptic curves~$E$
with even parity and conductor less than~$X$ and~$L(E,1)=0$.
Then $\tilde R(X)\sim cX^{19/24}(\log X)^{3/8}$
for some explicit constant~$c>0$.
\end{heuristic}

From Guess~\ref{GUESS} we get that the number of even parity curves
with $0<\Delta<qX$ and $D/N=q$ and $L(E,1)=0$ is given by
$$\hat W(qX)\cdot\bigl(\eta(q)/\beta_\tau\bigr)
\cdot{\rm Prob}\bigl[D/N=q\bigr],$$
and we sum this over all~$q$.
As we argued above, the tail of the sum does not affect the asymptotic
(and so we can take~$\log\Delta\sim\log N$ in~$\hat W$),
and again we get that the $q$-sum converges.
This then gives the desired asymptotic for the number of even parity
curves with conductor less than $X$ and vanishing central $L$-value
(upon arguing similarly for curves with negative discriminant).

\section{Torsion and isogenies\label{section:torsion}}
We can also count curves that have a given torsion group or isogeny
structure. For instance, an elliptic curve with a 2-torsion point
can be written as an integral model in the form $y^2=x^3+ax^2+bx$
where $\Delta=16b^2(a^2-4b)$; thus, by lattice-point counting,
we estimate about $\sqrt X$ curves
with absolute discriminant less than~$X$. The effect on the conductor
can perhaps more easily be seen by using the Fricke parametrisation
$c_4=(t+16)(t+64)T^2$ and $c_6=(t-8)(t+64)^2T^3$ of curves with
a rational \hbox{$2$-isogeny,} and then substituting $t=p/q$ and $V=T/q$
to get $c_4=(p+16q)(p+64q)V^2$ and $c_6=(p-8q)(p+64q)^2V^3$ so that
$\Delta=p(p+64q)^3q^2V^6$. The summation over the
twisting parameter~$V$ just multiplies our estimate by a constant,
while ABC-estimates imply that there should be no more than
$X^{2/3+\epsilon}$ coprime pairs $(p,q)$ with the squarefree kernel
of $pq(p+64q)$ smaller than $X$ in absolute value.
So we get the heuristic that almost all curves have no $2$-torsion, 
even under ordering by conductor. Indeed, the exceptional set is
so sparse that we can ignore it in our calculations.
A similar argument applies for other isogenies, and more
generally for splitting of division polynomials.
Also, the results of Duke~\cite{duke} for exceptional
primes are applicable here, albeit with a different ordering.

\section{Experiments}
We wish to provide some experimental data for the above heuristics.
In particular, the two large datasets of
Brumer-McGuinness \cite{brumer-mcguinness}
and Stein-Watkins \cite{stein-watkins}
for curves of prime conductor
up to $10^8$ and $10^{10}$ show little drop
in the observed average (analytic) rank.
Brumer and McGuinness considered about 310700 curves with prime conductor
less than $10^8$ and found an average rank of about $0.978$,
while Stein and Watkins extended this to over 11 million curves with
prime conductor up to $10^{10}$ and found an average rank of about~$0.964$.
Both datasets are expected to be nearly exhaustive\footnote
{This is one reason to take prime conductor curves;
 we also have $|\Delta|=N$ with few exceptions.}
amongst curves with prime conductor up to the given limit.
These results led some to speculate that
the average rank might be constant.
To test this, we chose a selection of curves with prime conductor
of size~$10^{14}$. It is non-trivial to get a good dataset, since
we must account for congruence conditions on the elliptic curve
coefficients and the variation of the size of the real period.

\subsection
{Average analytic rank for curves with prime conductor near $10^{14}$}
As in~\cite{stein-watkins},
we divided the $(c_4,c_6)$ pairs into 288 congruence classes with
$(\tilde c_4,\tilde c_6)=
\bigl(c_4\>\text{mod}\>576,c_6\>\text{mod}\>1728\bigr)$.
Many of these classes force the prime 2 to divide the discriminant,
and thus do not produce any curves of prime conductor.
For each class $(\tilde c_4,\tilde c_6)$, we took the
10000 parameter selections
$$(c_4,c_6)=\bigl(576(1000+i)+\tilde c_4,1728(100000+j)+\tilde c_6\bigr)
\>\text{for}\> (i,j)\in [1..10]\times [1..1000],$$
and then of these 2880000 curves,
took the 89913 models that had prime discriminant
(note that all the discriminants are positive).
This gives us good distribution across congruence classes,
and while the real period does not vary as much as possible,
below we will attempt to understand how this affects the average rank.

It then took a few months to compute the (suspected) analytic ranks for
these curves. We got about $0.937$ for the average rank. We then
did a similar experiment for curves with negative discriminant given by
$$(c_4,c_6)=\bigl(576(-883+i)+\tilde c_4,1728(100000+j)+\tilde c_6\bigr)
\>\text{for}\>(i,j)\in [1..10]\times [1..1000],$$
took the subset of 89749 curves with prime conductor,
and found the average rank to be about~$0.869$.
This discrepancy between positive and negative discriminant is
also in the Brumer-McGuinness and Stein-Watkins datasets, and
indeed was noted in \cite{brumer-mcguinness}.\footnote
{``An interesting phenomenon was the systematic influence of the
 discriminant sign on all aspects of the arithmetic of the curve.''}
We do not average the results from positive and negative
discriminant; the Brumer-McGuinness Conjecture \ref{conj:bmcg}
implies that the split is not~50-50.

In any case, our results show a substantial drop in the average rank,
which, at the very least, indicates that the average rank is not constant.
The alternative statistic of frequency of positive rank for curves with
even parity also showed a significant drop. For prime positive discriminant
curves it was 44.1\% for Brumer-McGuinness and 41.7\% for Stein-Watkins,
but only 36.0\% for our dataset --- for negative discriminant curves,
these numbers are 37.7\%, 36.4\%, and 31.3\%.

\subsection{Variation of real period}
Our random sampling of curves with prime conductor of size $10^{14}$
must account for various properties of the curves if our results
are to possess legitimacy. Above we speculated that
the real period plays the most significant r\^ole,
and so we wish to understand how our choice has affected it.

To judge the effect that variation of the real period might have,
we did some comparisons with the Stein-Watkins database.
First consider curves of positive prime discriminant,
and write $E$ as $y^2=4x^3+b_2x^2+2b_4x+b_6$
and $e_1>e_2>e_3$ for the real roots of the cubic.
We looked at curves with even parity and
considered the frequency of positive rank as a function
of the root quotient $t={e_1-e_2\over e_1-e_3}$, noting that\footnote
{The calculation follows as in the previous sections;
 via calculus, we can compute that this function is maximised
 at $t\approx .0388505246188$ with a maximum just below $4.414499094$.}
$\Omega_{\rm re}\Delta^{1/12}={2^{1/3}\pi(t-t^2)^{1/6}\over\agm(1,\sqrt t)}$.
The curves we considered all had $0.617<t<0.629$.

\begin{figure}[h]
\begin{center}
\scalebox{0.97}{\includegraphics{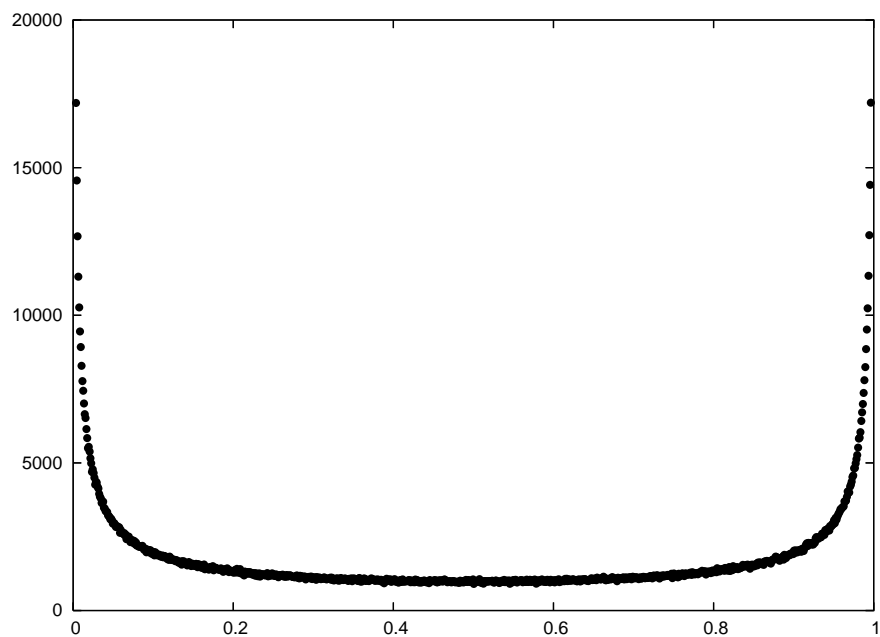}}
\end{center}
\vskip-12pt\noindent
\caption{$\Delta>0$: Curve distribution as a function of~$t$
\label{fig:posdisc-rootquotient}}
\end{figure}

However, similar to when we considered curves ordered by conductor,
before counting curves with extra rank,
we should first simply count curves. Figure~\ref{fig:posdisc-rootquotient}
indicates the distribution of root quotient $t$ for the
curves of prime (positive) discriminant and even parity
from the Stein-Watkins database (more than 2 million curves
meet the criteria). The $x$-axis is divided up into bins
of size $1/1000$; there are more than 100 times as many curves with $t<0.001$
as there with $0.500<t<0.501$, with the most extremal dots not
even appearing on the graph.

\vskip0pt\noindent
\vspace*{-1ex}
\begin{figure}[H]
\begin{center}
\scalebox{0.97}{\includegraphics{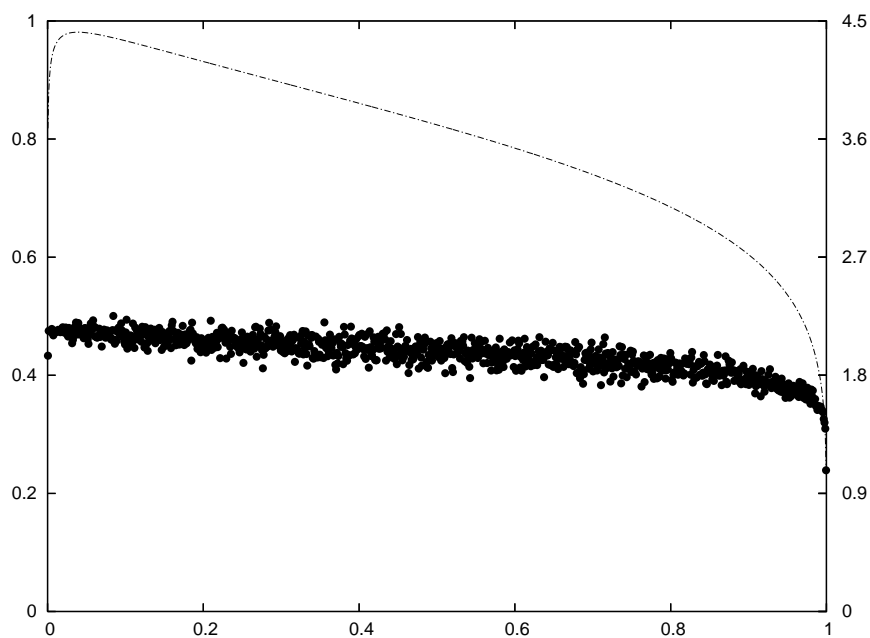}}
\end{center}
\vskip-12pt\noindent
\caption
{$\Delta>0$:\kern-0.5pt{} {}Positive rank frequency as a function of the
\newline\mbox{}\hspace{83pt}
root quotient~$t$, and $\Omega_{\rm re}\Delta^{1/12}$ as a function of~$t$.
\label{fig:posdisc-rankfreq}}
\end{figure}

Next we plot the frequency of $L(E,1)=0$ as a function of the root quotient
in Figure~\ref{fig:posdisc-rankfreq}.
Since there are only about 1000 curves in some of our bins, we do not
get such a nice graph. Note that the left-most and especially the
rightmost dots are much below their nearest neighbors,
the graph slopes down in general, and drops more at the end.
We see no evidence that our results should be overly biased.
In particular, the frequency of $L(E,1)=0$
is 41.7\% amongst all even parity curves of
prime discriminant in the Stein-Watkins database,
and is 42.8\% for the 12324 such curves with $0.617<t<0.629$.
The function plotted (labelled on the right axis)
in Figure~\ref{fig:posdisc-rankfreq} is of
$\Omega_{\rm re}\Delta^{1/12}={2^{1/3}\pi(t-t^2)^{1/6}\over\agm(1,\sqrt t)}$
as a function of~$t$, and note that this goes to 0 as $t \rightarrow 0,1$;
there is nothing canonical about the choice of our $t$-parameter,
and we chose it more for convenience than anything else.

\begin{figure}[H]
\begin{center}
\scalebox{0.97}{\includegraphics{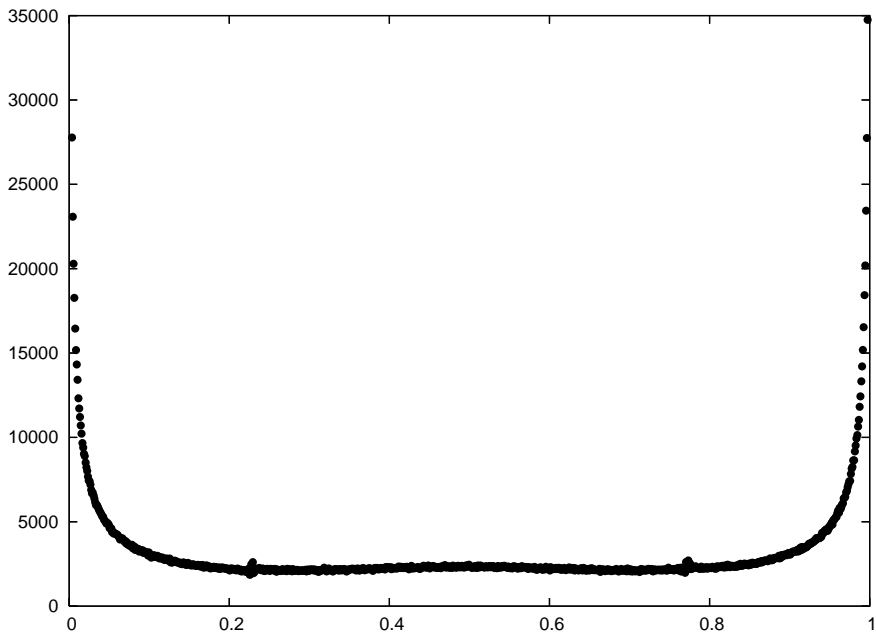}}
\end{center}
\vskip-18pt\noindent
\caption{$\Delta<0$: Distribution of curves as a function of $C$
\label{fig:negdisc-rootquotient}}
\vskip-6pt
\end{figure}

Similar computations can be made in the case of negative discriminant,
which we briefly discuss for completeness (again restricting to curves
with even parity where appropriate). Let $r$ be the real root of
the cubic polynomial $4x^3+b_2x^2+2b_4x+b_6$,
and $Z>0$ the imaginary part of the conjugate pair of nonreal roots.
Letting $\tilde r=r+b_2/12$ and $c=\tilde r/Z$ we then have\footnote
{This is maximised at $c\approx -33.58515148525$,
 with the maximum a bit less than $8.82921518$.}
\vskip-2pt\noindent
$$\Omega_{\rm re}|\Delta|^{1/12}=
{\pi\sqrt{2}\over (1+9c^2/4)^{1/12}
\agm\Bigl(1,\sqrt{{1\over 2}+{3c\over 4\sqrt{1+9c^2/4}}}\Bigr)}.$$
We renormalise via taking~$C=1/2+\arctan(c)/\pi$,
and graph the distribution of curves versus~$C$ in
Figure~\ref{fig:negdisc-rootquotient}.
The symmetry\footnote
{The blotches around 0.22-0.23 and 0.77-0.78 appear to come from
 the fact that curves with $a_4$ small (in particular~$\pm 1$)
 tend to have $C$ in these ranges (for our discriminant range),
 and this causes instability in the counting function.}
of the graph might indicate that the coordinate transform is reasonable.
All our curves have~$0.555<C<0.557$.

\vskip-8pt\noindent

\begin{figure}[h]
\begin{center}
\scalebox{0.97}{\includegraphics{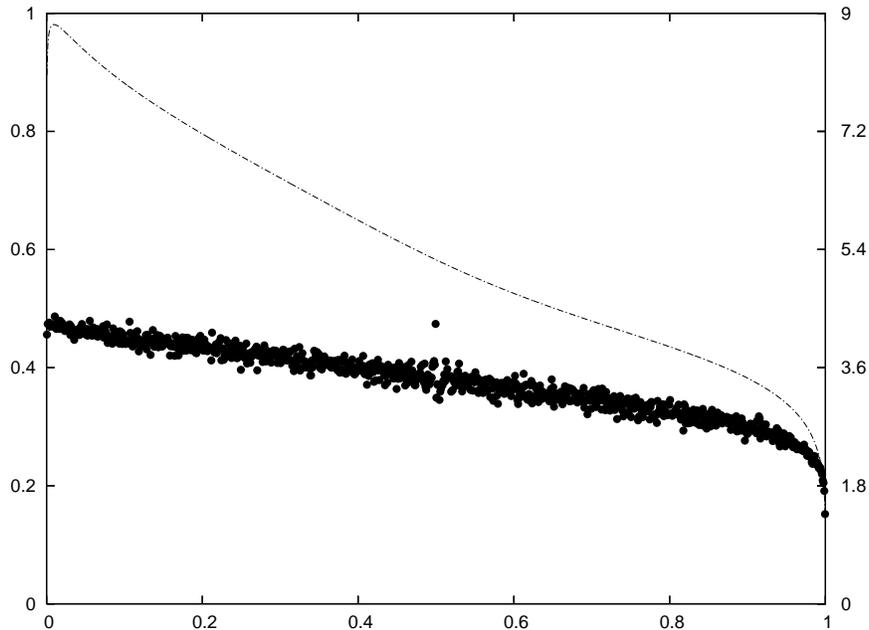}}
\end{center}
\vskip-12pt\noindent
\caption{$\Delta<0$: Positive rank frequency as a function of~$C$,
\newline\mbox{}\hspace{82.5pt}
and $\Omega_{\rm re}|\Delta|^{1/12}$ as a function of~$C$.
\label{fig:negdisc-rankfreq}}
\vskip-6pt\noindent
\end{figure}

Next we plot the frequency of $L(E,1)=0$ as a function of the root quotient
in Figure~\ref{fig:negdisc-rankfreq}. Again we also graph the function
$\Omega_{\rm re}|\Delta|^{1/12}$ on the right axis.
Here the drop-off is more pronounced
than with the curves of positive discriminant.
Note the floating dot around $C=1/2$.
Indeed the 100 closest curves with $C<1/2$ all have positive rank;
this breaks down when crossing the $1/2$-barrier.
This is not particularly a mystery; these curves have $a_6=0$ and/or $b_6=1$,
and thus have an obvious rational point.
Recall that $C=1/2$ corresponds to $c=0=\tilde r$.

We again see no evidence that our results should be biased.
In particular, the frequency of $L(E,1)=0$ is 36.4\% amongst
all even parity curves of negative prime discriminant
in the Stein-Watkins database, and is 37.0\% for the
4695 such curves with $0.555<C<0.557$.

\subsection{Other considerations}
The idea that the ``probability'' that an even parity curve possesses
positive rank should be proportional to $\sqrt{\Omega_{\rm re}}$
is perhaps overly simplistic --- in particular, it is not borne out
too precisely by the Stein-Watkins dataset.
We consider positive prime discriminant curves with even parity;
for those with $0.64<\Omega_{\rm re}<0.65$ we have 78784 curves of
which 45.9\% have positive rank, while of the 9872 with
$0.32<\Omega_{\rm re}<0.325$ we have 36.0\% with positive rank,
for a ratio of 1.28, which is not too close to~$\sqrt 2$.
One consideration here is that we have placed a discriminant limit
on our curves, and there are curves with larger discriminant and
$0.32<\Omega_{\rm re}<0.325$ that we have not considered.
This, however, is extra-particular to the idea that only the
real period should be of import.

One possibility is that curves with small discriminant and/or
large real period have \textit{smaller} probability of $L(E,1)=0$
that our estimate of $c\sqrt{\Omega_{\rm re}}$ would suggest --- indeed,
it might be argued (maybe due to arithmetic considerations,
or perhaps explicit formulae for the zeros of $L$-functions)
that curves with such small discriminant cannot realise their
nominal expected frequency of positive rank. Unfortunately,
we cannot do much to quantify these musings, as the effect would
likely be in a secondary term, making it difficult to detect experimentally.
Note also that a relative depression of rank for small discriminant
curves would give a reason for the near-constant average rank
observed by Brumer-McGuinness and Stein-Watkins.

\subsection{Mordell-Weil lattice distribution for rank 2 curves}
We have other evidence that curves of small discriminant might not
behave quite as expected. We undertook to compute generators for
the Mordell-Weil group for all 2143079 curves of (analytic) rank 2
of prime conductor less than $10^{10}$ in the Stein-Watkins
database.\footnote
{We also computed the Mordell-Weil group for curves with higher ranks
 but do not describe the obtained data here.}
J. E. Cremona ran his \verb+mwrank+ programme \cite{mwrank} 
on all these curves, and it was successful in provably finding
the Mordell-Weil group for 2114188 of these.
For about 2500 curves, the search region was
too big to find the 2-covering quartics via invariant methods,
while around 8500 curves had a generator of large height that
could not be found, and over 18000 had 2-Selmer rank greater than~2.
We then used the \verb+FourDescent+ machinery of MAGMA \cite{magma}
which reduced the number of problematic curves to 54. Of these,
19 have analytic \Sha\ of 16.0 and we expect that either 3-descent
or 8-descent \cite{stamminger} will complete (assuming GRH to
compute the class group) the Mordell-Weil group verification;
for the 35 other curves, there is likely a generator
of height more than 225 which we did not attempt to find.\footnote
{A bit more searching might resolve a few of the outstanding cases,
 but the extremal case of $[0,0,1,-237882589,-1412186639384]$ appears
 to have a generator of height more than 600, and thus other methods
 will likely be needed to try to find it. T.~A. Fisher has recently
 used 6-descent to find some of the missing points.}

We then looked at the distribution of the Mordell-Weil lattices
obtained from the induced inner product from the height pairing;
since all of our curves have rank 2, we get 2-dimensional lattices.
We are not so interested in the size of the obtained lattices, but
more so in their shape. Via the use of lattice reduction (which
reduces to continued fractions in this case), given any two generators
we can find the point $P$ of smallest positive height on the curve.
By normalising $P$ to be the unit vector, we then get a vector in the
upper-half-plane corresponding to another generator~$Q$. Via the standard
reduction algorithm, we can translate $Q$ so that it corresponds to a 
point in the fundamental domain for the action of~${\rm SL}_2({\bf Z})$.
Finally, by replacing $Q$ by $-Q$ if necessary, we can ensure that this
point is in the right half of the fundamental domain
(in other words, we must choose an embedding for our Mordell-Weil lattice).
In this manner, for each rank 2 curve we associate a unique point
$z=x+iy$ in the upper-half-plane with $x^2+y^2\ge 1$ and $0\le x\le 1/2$.

With no other guidance, we might expect that the obtained distribution
for the $z$ is given by\footnote
{Siegel \cite{siegel} similarly uses Haar measure to put a natural
 measure on $n$-dimensional lattices of determinant~$1$.}
the Haar measure $(dx\, dy)/y^2$.
We find, however, that this is not borne out too well by experiment.
In particular, we should expect that
${1/2\over \pi/6}\approx 95.5\%$ of the curves should have $y\ge 1$,
while the experimental result is about~93.5\%. Furthermore, we should
expect that the proportion of curves with $y\ge Y$ should die off like
$1/Y$ as $y\rightarrow\infty$; however, we get that 35.4\% of the curves
have $y\ge 2$, only 9.4\% have $y\ge 4$, while 1.7\% have $y\ge 8$
and 0.2\% have $y\ge 16$. The validity of the vertical distribution
data might be arguable based upon concerns regarding the discriminant
cutoff of our dataset, but the horizontal distribution is also skewed.
If we consider only curves with~$y\ge 1$, then we should get uniform
distribution in the $x$-aspect; however, Table~\ref{table:xdist} shows that
we do not have such uniformity.

\begin{table}[h]
\caption{Horiztonal distribution of rank 2 lattices with $y\ge 1$
\label{table:xdist}}
\begin{center}
\begin{tabular}{*{4}{|@{\hskip6pt}c@{\hskip6pt}}|}\hline
$0.00\le x<0.05$&9.0\%&$0.25\le x<0.30$&10.0\%\\
$0.05\le x<0.10$&9.6\%&$0.30\le x<0.35$&10.2\%\\
$0.10\le x<0.15$&9.8\%&$0.35\le x<0.40$&10.5\%\\
$0.15\le x<0.20$&9.9\%&$0.40\le x<0.45$&10.6\%\\
$0.20\le x<0.25$&10.0\%&$0.45\le x\le 0.50$&10.5\%\\\hline
\end{tabular}
\end{center}
\end{table}

We cannot say whether these unexpected results from the experimental data
are artifacts of choosing curves with small discriminant; it is just as
probable that our Haar-measure hypothesis concerning the lattice distribution
is simply incorrect.

\subsection{Symmetric power L-functions}
Similar to questions about the vanishing of $L(E,s)$, we can ask about
the vanishing of the symmetric power $L$-functions $L({\rm Sym}^{2k-1} E,s)$.
We refer the reader to \cite{martin-watkins} for more details about this,
but mention that, due to conjectures of Deligne and more generally
Bloch and Be\u\i linson \cite{beilinson},
we expect that we should have a formula similar
to that of Birch and Swinnerton-Dyer, stating that
$L({\rm Sym}^{2k-1} E,k)(2\pi N)^{k\choose 2}/
\Omega_+^{k+1\choose 2}\Omega_-^{k\choose 2}$
should be rational with small denominator. Here, for $k$ odd,
$\Omega_+$ is the real period and $\Omega_-$ the imaginary period,
with this reversed for $k$ even. Ignoring the contribution
from the conductor, and crudely estimating that
$\Omega_+\approx\Omega_-\approx 1/\Delta^{1/12}$, an application
of discretisation as before gives that the probability
that $L({\rm Sym}^{2k-1} E,s)$ has even parity and
$L({\rm Sym}^{2k-1} E,k)=0$ is bounded above (cf.~the ignoring of~$N$)
by $c(\log\Delta)^{3/8}\cdot\sqrt{1/\Delta^{k^2/12}}$.
Again following the analogy of above, we can then upper-bound
the number of curves with conductor less than~$X$
with even-signed symmetric $(2k-1)$st power and $L({\rm Sym}^{2k-1} E,k)=0$
by $c_k(\epsilon) X^{5/6-k^2/24+\epsilon}$ for every $\epsilon>0$.
It could be argued that
we should order curves according to the conductor of the symmetric
power $L$-function rather than that of the curve, but we do not think
such concerns are that relevant to our imprecise discussion.
In particular, the above estimate predicts that there are
finitely many curves with extra vanishing when~$k\ge 5$. It should be
said that this heuristic will likely mislead us about curves with
complex multiplication, for which the symmetric power $L$-function
factors (it is imprimitive in the sense of the Selberg class),
with each factor having a 50\% chance of having odd parity. However,
even ignoring CM curves, the data of \cite{martin-watkins} find a handful
of curves for which the 9th, 11th and even the 13th symmetric powers appear
(to 12 digits of precision)
to have a central zero of order~2. We find this surprising, and casts some
doubt about the validity of our methodology of modelling of vanishings.

\subsection{Quadratic twists of higher symmetric powers}
The techniques we used earlier in this paper have also
been used to model vanishings in quadratic twist families,
and we can extend the analyses to symmetric powers.

\subsubsection{Non-CM curves}
We fix a non-CM curve~$E$ and let $E_d$ be its $d$th quadratic twist,
taking $d$ to be a fundamental discriminant. From an analogue of the
Birch--Swinnerton-Dyer conjecture we expect to get a small-denominator
rational from the quotient\footnote
{The contribution from the conductor actually comes from non-integral
 Tamagawa numbers from the Bloch-Kato exponential map, and in the case
 of quadratic twists, the twisting parameter~$d$ should not appear in
 the final expression.} 
$L({\rm Sym}^3 E_d,2)(2\pi N_E)/\Omega_{\rm im}(E_d)^3\Omega_{\rm re}(E_d)$.
We have that $\Omega_{\rm im}(E_d)^3\Omega_{\rm re}(E_d)
\approx \Omega_{\rm im}(E)/d^{3/2}\cdot \Omega_{\rm re}(E)/d^{1/2}$
and so we expect the number of fundamental discriminants $|d|<D$
such that $L({\rm Sym}^3 E_d,s)$ has even parity with
$L({\rm Sym}^3 E_d,2)=0$ to be given crudely (up to log-factors)
by $\sum_{d<D} c/\sqrt{d^2}$.
So we expect about $(\log D)^b$ quadratic twists with double zeros
for the 3rd symmetric power; generalising predicts finitely many
extra vanishings for higher (odd) powers.

\begin{table}[h]
\caption
{Fundamental $d$ with
 $\mathop{\rm ord}\limits_{s=2} L({\rm Sym}^3 E_d,s)\ge 2$\label{tbl:Ed}}
\vspace*{-12pt}
\begin{center}
\begin{tabular}{|r|l|}\hline
11a& $-40$ $-52$ $-563$ $-824$ $-1007$
$-1239$ $-1460$ $-1668$ $-1799$ $-2207$\\
&$-2595$ $-2724$ $-2980$ $-3108$ $-3592$ $-4164$ $-4215$ $-4351$ $-4399$\\
&12 69 152 181 232 273 364 401 412 421 444 476 488 652 669 696\\
&933 1101 1149 1401 1576 1596 1676 1884 1928 2348 2445 2616\\
&2632 3228 3293 3404 $3720^\star$ 3793 4060 4093 4161 4481 4665 4953\\\hline
14a&
$-31$ $-52$ $-67$ $-87$ $-91$ $-111$ $-203$
$-223$ $-255$ $-264$ $-271$ $-311$ $-327$\\
&$-367$ $-535$ $-552$ $-651$ $-759$ $-804$
$-831$ $-851$ $-852$ $-920$ $-1099$ $-1263$\\
&$-1267$ $-1335$ $-1524$ $-1547$ $-1567$
 $-1623$ $-1679$ $-1707^\star$ $-2047$ $-2235$\\
&$-2280$ $-2407$ $-2443$ $-2563$ $-2824$
 $-2831$ $-3127$ $-3135$ $-3523$ $-4119$\\
&$-4179$ $-4191$ $-4323$\\
&137 $141^\star$ 229 233 281 345 469 473
 492 $497^\star$ 697 901 1065 1068 1353 1457\\
&1481 1513 1537 1793 1873 1905 2024 2093 2193 2265 2321 2589 2657 2668\\
&2732 2921 2981 2993 3437 3473 3529 4001 4124 4389 4488 4661 4817\\\hline
15a&
$-11$ $-51$ $-71$ $-164$ $-219$ $-232$ $-292$ $-295$ $-323^\star$ $-340$
$-356$ $-399$ $-519$\\
&$-580$ $-583$ $-584$ $-671$ $-763$ $-804$ $-851$
$-879$ $-943$ $-1012$ $-1060$ $-1151$\\
&$-1199$ $-1284$ $-1288$ $-1363$
$-1551$ $-1615$ $-1723$ $-1732$ $-2279$ $-2291$\\
&$-2379$ $-2395$
$-2407$ $-2571$ $-2632$ $-2635$ $-2756$ $-3396$ $-3588^\star$ $-3832$\\
&17 21 61 77 136 156 181 229 349 444 481 501 545 589 649 781 876 905\\
&924 949 1009 1144 1249 1441 1501 1580 1621 1804 1861 1921 2041 2089\\
&2109 2329 2581 2829 2840 2933 3001 3916\\\hline
\end{tabular}
\end{center}
\end{table}

We took the curves 11a:$[0,-1,1,0,0]$ and 14a:$[1,0,1,-1,0]$,
and computed either $L({\rm Sym}^3 E_d,2)$ or $L'({\rm Sym}^3 E_d,2)$
for all fundamental discriminant $d$ with~$|d|<5000$.
We did the same for 15a:$[1,1,1,0,0]$ for~$|d|<4000$.
We then looked at the number of vanishings (to 9 digits of precision).
For 11a we found 58 double zeros and one triple zero
(indicated by a star in Table~\ref{tbl:Ed}) while for 14a
we found 88 double zeros and three triple zeros, and 15a yielded
83 double zeros and two triple zeros.

\subsubsection{CM curves}
Next we consider CM curves, for which we can compute significantly more data,
but the modelling of vanishings is slightly different.
Let $E$ be a rational elliptic curve with~CM, and $\psi$ its Hecke character.
We shall take~$\psi$ to be ``twist-minimal'' --- this is not the same
as the ``canonical'' character of Rohrlich~\cite{rohrlich},
but rather we just take $E$ to be a minimal (quadratic) twist.
Indeed, we shall only consider 11 different choices of~$E$,
given (up to isogeny class) by 27a, 32a, 36a, 49a, 121a,
256a, 256b, 361a, 1849a, 4489a, and~26569a, noting that 27a/36a
and 32a/256b are respectively cubic and quartic twist-pairs.
In our tables, these can appear in a briefer format,
such as $67^2$ for~4489a.

We normalise the Hecke $L$-function $L(\psi,s)$ to have $s=1$
be the center of the critical strip.
For $d$ a fundamental discriminant, we let $\psi_d$
be the Hecke Gr\"ossencharacter $\psi$ twisted by the quadratic
Dirichlet character of conductor~$d$. Finally, note that
the symmetric powers $L({\rm Sym}^k \psi,s)$ are just~$L(\psi^k,s)$,
where we must take $\psi^k$ to be the primitive underlying Gr\"ossencharacter.

We then expect $L(\psi^3,2)(2\pi)/\Omega_{\rm im}(E)^3$
to be rational with small denominator.
We can then use discretisation as before to count the expected number
of fundamental discriminants $|d|<D$ for which the
$L$-function $L(\psi_d^3,s)$ has even parity
but vanishes at the central point --- since we have
$\Omega_{\rm im}(E_d)^3\approx \Omega_{\rm im}(E)/d^{3/2}$,
we expect the number of discriminants $d$ that
yield even parity and $L(\psi_d^3,2)=0$
is crudely given by $\sum_{d<D} 1/\sqrt{d^{3/2}}$,
so we should get about $D^{1/4}$ such discriminants up to~$D$.

For higher symmetric powers, we expect that
$L(\psi^{2k-1},k)(2\pi)^{k-1}/\Omega_+(E)^{2k-1}$
is rational with small denominator, and thus get that
there should be finitely many quadratic twists of even parity
with vanishing central value.

We took the above eleven CM curves and took their (fundamental) quadratic
twists up to~$10^5$. We must be careful to exclude twists that are
isogenous to other twists. In particular, we need to define a
{\bf primitive} discriminant for a curve with CM by an order
of the field $K$ --- this is a fundamental discriminant~$d$
such that ${\rm disc}(K)$ does not divide~$d$,
expect for $K={\bf Q}(i)$ when $d>0$ is additionally primitive when~$8\|d$.
Note also that 27a and 36a have the same symmetric cube $L$-function.

\vskip-12pt\noindent
\begin{table}[h]
\caption{Counts of double order zeros for primitive twists\label{tbl:overview}}
\vskip-8pt\noindent
\begin{center}
\begin{tabular}{|c|ccccccccccc|}\hline
&27a&32a&36a&49a&121a&256a&256b&361a&1849a&4489a&26569a\\\hline
3rd&59&32&-&67&78&32&21&45&28&31&1\\
5th&3&1&5&2&1&2&2&0&0&0&0\\
7th&0&0&2&0&1&0&0&0&0&0&0\\\hline
\end{tabular}
\end{center}
\end{table}
\vskip-6pt\noindent

Table~\ref{tbl:overview} lists the following results for counts
of central double zeros (to 32 digits) for the $L$-functions
of the 3rd, 5th, and 7th symmetric powers.\footnote
{We found no even twists with~$L(\psi_d^9,5)=0$,
 and no triple zeros appeared in the data.}
Tables~\ref{tbl:data} and \ref{tbl:higher}
list the primitive discriminants that yield
the double zeros. The notable signedness can be
explained via the sign of the functional equation.\footnote
{The local signs at $p=2,3$ involve wild ramification are more complicated
 (see~\cite{whitehouse,kobayashi,dm} for a theoretical description),
 and thus there is no complete correlation in some cases.}
We are unable to explain the paucity of double zeros for twists of~26569a;
Liu and Xu have the latest results~\cite{liu-xu} on the vanishing of
such $L$-functions, but their bounds are far from the observed data.
Similarly, the last-listed double zero for 4489a at 67260 seems quite small.

There appear to be implications vis-a-vis higher vanishings in some cases;
for instance, except for 27a, in the thirteen cases that $L(\psi_d^5,s)$ has
a double zero at $s=3$ then $L(\psi_d,s)$ also has a double zero at~$s=1$.
Similarly, the 7th symmetric power for the 27365th twist of 121a has a
double zero, as does the 3rd symmetric power, while the $L$-function
of the twist itself has a triple zero.
Also, the 22909th twist of 36a has double zeros
for its first, third, and fifth powers (note that 36a does not appear
in Table~\ref{tbl:data} as the data are identical to that for~27a).

\begin{table}[H]
\vspace*{-18pt}
\caption{Primitive $d$ with
         $\mathop{\rm ord\,}\limits_{s=2} L(\psi_d^3,s)=2$
\label{tbl:data}}
\vspace*{-6pt}
\begin{center}
\begin{tabular}{|r|l|}\hline
27a&172 524 1292 1564 1793 3016 4169 4648 6508 9149 9452 9560 10636\\
&11137 12040 13784 14284 15713 17485 17884 22841 22909 22936 25729\\
&27065 27628 29165 30392 34220 35749 38636 40108 41756 44221 47260\\
&51512 54385 57548 58933 58936 58984 59836 59996 62353 64268 70253\\
&74305 77320 77672 78572 84616 86609 86812 87013 92057 95861 96556\\
&97237 99817\\\hline
32a& $-395$ $-5115$ $-17803$ $-25987$ $-58123$\\
& $-60347$ $-73635$ $-79779$ $-84651$ $-99619$\\
& 257 1217 2201 2465 14585 26265 45201 82945\\
&4632 5336 5720 7480 9560 30328 30360\\
&31832 38936 45848 69784 71832 83512 92312\\\hline
49a&
$-79$ $-311$ $-319$ $-516$ $-856$ $-1007$
$-1039$ $-1243$ $-1391$ $-1507$ $-1795$\\
&$-2024$ $-2392$ $-2756$ $-2923$ $-3527$
$-3624$ $-4087$ $-4371$ $-4583$ $-4727$\\
&$-5431$ $-5524$ $-5627$ $-6740$ $-7167$
$-7871$ $-8095$ $-8283$ $-10391$ $-10628$\\
& $-13407$ $-13656$ $-13780$ $-16980$ $-18091$
$-22499$ $-27579$ $-28596$ $-30083$\\
& $-30616$ $-32303$ $-32615$ $-36311$
$-36399$ $-38643$ $-39127$ $-40127$ $-42324$\\
& $-52863$ $-64031$ $-64399$ $-66091$ $-66776$
 $-66967$ $-69647$ $-70376$ $-71455$\\
& $-72663$ $-73487$ $-73559$ $-77039$ $-84383$
$-90667$ $-91171$ $-98655$ $-98927$\\\hline
$11^2$&
12 140 632 1160 1208 1308 1704 1884 2072 2136 2380 2693 2716 3045\\
&4120 4121 5052 5528 5673 5820 6572 7532 11053 11208 12277 12568\\
&12949 13884 14844 15465 16136 18588 18885 19020 19884 24060 25788\\
&27365 27597 28265 28668 29109 29573 32808 32828 35261 36552 37164\\
&38121 38297 44232 44873 49512 49765 50945 52392 54732 55708 56076\\
&56721 58460 59340 65564 66072 66833 71688 72968 79557 80040 80184\\
&83388 84504 84620 84945 86997 87576 92460 95241\\\hline
\kern-1.5pt 256a&
401 497 2513 3036 3813 6933 6941 9596 9932 11436 14721 17133 17309\\
&18469 21345 21749 26381 26933 28993 29973 30461 33740 51469 53084\\
&62556 63980 67721 69513 73868 76241 81164 87697\\\hline
\kern-1.5pt 256b&
73 345 3521 5133 6693 7293 21752 25437 27113 34657 38485 41656\\
&42433 44088 46045 75581 79205 83480 89737 93624 96193\\\hline
$19^2$&
44 60 1429 1793 3297 3340 3532 3837 3880 4109 5228 5628 7761 8808\\
&9080 9388 12280 12313 12545 13373 13516 13897 19164 22204 23241\\
&25036 25653 41205 41480 42665 43429 44121 44285 44508 45660 48828\\
&50584 52989 64037 74585 75324 76921 81885 85036 96220\\\hline
$43^2$&
88 152 440 2044 4268 5852 6376 7880 8908 9880 14252\\
&15681 17864 20085 20353 28492 29477 45368 55948 56172\\
&57409 60177 68136 79916 84524 85580 86853 96216\\\hline
$67^2$&
17 57 869 1612 1628 3260 6380 6385 7469 8328 11017 13772\\
&14152 14268 14552 15901 22513 24605 24664 27992 29676 33541\\
&33789 36344 36588 38028 40280 43041 49884 62353 67260 \\\hline
\kern-1.5pt$163^2$&30720\\\hline
\end{tabular}
\end{center}
\vspace*{-7.5pt}
\caption{Primitive $d$ with
       $\mathop{\rm ord\,}\limits_{s=k} L(\psi_d^{2k-1},s)=2$
         for some $k\ge 3$\label{tbl:higher}}
\vspace*{-2pt}
\begin{center}
\begin{tabular}{|r|l|r|l|}\hline
27a&5th: $-13091$ 4040 18044&49a&5th: 437 19317\\
32a&5th: 1704&121a&5th: $-183$\quad 7th: 27365\\
36a&5th: $-856$ $-2104$ $-31592$ $-88580$ 22909&256a&5th: $-79$ $-21252$\\
36a&7th: $-95$ 2488&256b&5th: $-511$ 89320\\\hline
\end{tabular}
\end{center}
\end{table}
\goodbreak

\subsubsection{Comparison between the CM and non-CM cases}
For the twist computations for the symmetric powers,
we can go much further (about 20 times as far) in the CM case
because the conductors do not grow as rapidly.\footnote
{In \cite[\S 8]{RVZ} Rodriguez Villegas and Zagier
 mention the possibility of a Waldspurger-type formula
 for the twists of the Hecke Gr\"ossencharacters,
 but it does not seem that such a formula has ever appeared.
 Similarly, one might hope to extend the work of Coates and Wiles~\cite{cw}
 and/or Gross and Zagier~\cite{gz} to powers of Gr\"ossencharacters;
 there is some early work (among others) of Damerell \cite{damerell}
 in this regard, while Guo \cite{guo} shows partial results toward the
 Bloch-Kato conjecture.}
For the 3rd symmetric power, the crude prediction is that we should have
(asymptotically) many more extra vanishings for twists in the CM case
than in the non-CM case, but this is not borne out by the data.
Additionally, we have no triple zeros in the CM case (where the
dataset is almost 100 times as large),
while we already have six for the non-CM curves.
This is directly antithetical to our suspicion that there should
be more extra vanishings in the CM case.
As before, this might cast some doubt on our methodology of modelling
of vanishings.

\section{Acknowledgements}
The author was partially supported by
Engineering and Physical Sciences Research Council (EPSRC)
grant GR/T00658/01 (United Kingdom).
He thanks N.~D.~Elkies, H.~A.~Helfgott, and A.~Venkatesh for useful comments,
and N.~P.~Jones for the reference~\cite{duke}.

\end{document}